\documentclass[secnum,leqno]{rims-bessatsu}
\usepackage{amssymb, amsmath}
\usepackage{mathrsfs}
\usepackage{pb-diagram}
\usepackage[all]{xy}
\usepackage{graphics}
\usepackage[dvips]{graphicx}
\usepackage{eepic}
\usepackage{epic}
\usepackage{enumerate,mathtools}
\usepackage[hyperindex=true,colorlinks=false]{hyperref}
\newtheorem{thm}{Theorem}[section]
\newtheorem{crl}[thm]{Corollary}

\newtheorem{lmm}[thm]{Lemma}

\newtheorem{prp}[thm]{Proposition}

\theoremstyle{definition}
\newtheorem{dfn}[thm]{Definition}
\newtheorem{exa}[thm]{Example}
\newtheorem{nota}[thm]{Notation}
\theoremstyle{remark}
\newtheorem{rem}[thm]{Remark}
\newenvironment{eq-text}
{\begin{equation} \begin{minipage}[t]{0.85\linewidth}}
{\end{minipage} \end{equation} \ignorespacesafterend}

\DeclareMathOperator{\N}{\mathbb{N}}

\DeclareMathOperator{\R}{\mathbb{R}}
\DeclareMathOperator{\C}{\mathbb{C}}

\newcommand{\ti}{\tilde}
\newcommand{\ii}{^{-1}}
\newcommand{\pa}{\partial}
\newcommand{\ee}{\mathrm e}
\newcommand{\defeq}{\coloneqq} 
\newcommand{\col}{\colon\thinspace}          
\newcommand{\ens}{\enspace}
\newcommand{\ie}{{\emph{i.e.}}\ }
\newcommand{\dd}{{\mathrm d}}
\newcommand{\un}[1]{{\underline{#1}}}
\newcommand{\ov}{\overline}

\newcommand{\pp}[1]{^{(#1)}}

\DeclarePairedDelimiter\floor{\lfloor}{\rfloor}%
\DeclarePairedDelimiter\abs{\lvert}{\rvert}%
\newcommand{\be}{\beta}
\newcommand{\eps}{\varepsilon}
\newcommand{\ph}{\varphi}
\newcommand{\Om}{\Omega}
\newcommand{\om}{\omega}
\newcommand{\Ga}{\Gamma}
\newcommand{\ga}{\gamma}
\newcommand{\De}{\Delta}
\newcommand{\de}{\delta}
\newcommand{\la}{\lambda}
\newcommand{\La}{\Lambda}
\newcommand{\ze}{\zeta}

\newcommand{\gD}{{\vec{\mathscr D}}}
\newcommand{\gE}{\mathscr E}

\newcommand{\gL}{\mathscr L}
\newcommand{\gR}{\mathscr R}
\newcommand{\gO}{\mathscr O}

\newcommand{\cB}{\mathcal{B}}
\newcommand{\cM}{\mathcal{M}}
\newcommand{\cN}{\mathcal{N}}
\newcommand{\cS}{\mathcal{S}}

\newcommand{\dfs}{d.f.s.}
\newcommand{\fp}{\mathfrak{p}}
\newcommand{\fq}{\mathfrak{q}}
\newcommand{\uO}{\un{0}}
\newcommand{\uU}{\un{U}}
\newcommand{\uga}{\un{\ga}}
\newcommand{\vuze}{\vec{\un\ze}}

\newcommand{\Rp}{\R_{\geq0}}
\newcommand{\RCV}{\operatorname{RCV}}
\newcommand{\dist}{\operatorname{dist}}

\newcommand{\Det}{\operatorname{det}}

\newcommand{\begla}{\begin{equation}}
\newcommand{\beglab}[1]{\begin{equation}	\label{#1}}
\newcommand{\edla}{\end{equation}}

\newcommand{\ID}{\mathop{\hbox{{\rm Id}}}\nolimits}

\newcommand{\isom}{\xrightarrow{\smash{\ensuremath{\sim}}}}

\newcommand\restr[2]{{
  \left.\kern-\nulldelimiterspace 
  #1 
  \vphantom{\big|} 
  \right|_{#2} 
  }}

\title{Nonlinear analysis with endlessly continuable functions}
\TitleHead{}          
\author{Shingo \textsc{Kamimoto}
        and 
David \textsc{Sauzin}\footnote{Laboratorio Fibonacci, CNRS--CRM Ennio De
  Giorgi SNS Pisa, Italy.
}}

\AuthorHead{Shingo Kamimoto and David Sauzin }         
\support{Supported by \underline{Grant-in-Aid for JSPS Fellows Grant Number 15J06019}, and
French National Research Agency \underline{reference ANR-12-BS01-0017}}  
\begin{document}
%



\thispagestyle{empty}

\begin{center}
\resizebox{\linewidth}{!}{
\textbf{\Large Nonlinear analysis with endlessly continuable functions}}
\end{center}

\bigskip

\begin{center}
{\large Shingo \textsc{Kamimoto} and David \textsc{Sauzin}}
\end{center}

\begin{abstract}

  We give estimates for the convolution product of an arbitrary number
  of endlessly continuable functions. This allows us to deal with
  nonlinear operations for the corresponding resurgent series,
  e.g. substitution into a convergent power series.

\end{abstract}

\bigskip


\section{Introduction}\label{sec:0}

This is an announcement of our forthcoming paper \cite{KS}, the main
subject of which is the ring structure of the space~$\gR$ of \emph{resurgent}
formal series.

Recall that~$\gR$ is the subspace of the ring of formal series $\C[[z]]$
defined as
\beglab{eqdefgR}
\gR \defeq \cB\ii(\C\de \oplus \hat\gR),
\edla
where
$\cB \col \C[[z]] \to \C\de \oplus \C[[\ze]]$
is the formal Borel transform, defined by
\beglab{eqdefcB}
\ti\ph(z) = \sum_{j=0}^{\infty} \ph_j z^{j}
\mapsto 
\cB(\ti\ph)(\ze) = \ph_0\de + \hat\ph(\ze),
\ens \text{with}\; \;
\hat\ph(\ze) \defeq \sum_{j=1}^{\infty} \ph_j \frac{\ze^{j-1}}{(j-1)!},
\edla
and $\hat\gR$ is the subspace of~$\C[[\ze]]$ consisting of
all convergent power series which are ``endlessly continuable''
in the following sense:

%
\begin{dfn}   \label{defendlesscont}
  A convergent power series $\hat\ph\in\C\{\ze\}$ is said to be
  \emph{endlessly continuable} if, for every real $L>0$, there exists a
  finite subset~$F_L$ of~$\C$ such that~$\hat\ph$ can be
  analytically continued along every Lipschitz path 
  $\ga \col [0,1] \to \C$
of length $<L$ 
such that $\ga(0)=0$ and $\ga\big( (0,1] \big) \subset \C\setminus F_L$.
\end{dfn}
%

Resurgence theory was invented by J.~\'Ecalle in the early 1980s
\cite{E} and has many applications in the study of holomorphic
dynamical systems, analytic differential equations, WKB analysis, etc.
Here, we are dealing with a simplified version, inspired by
\cite{CNP}, which is sufficient for most applications 
(our definition of endlessly continuable functions is slightly more
restrictive than the one in \cite{CNP}, which is itself less general
than the definition of ``functions without a cut'' on which is based
\cite{E}
-- see also \cite{DO}).


We are interested in nonlinear operations in the space of formal
series, like substitution of one or several series without constant
term $\ti\ph_1,\ldots, \ti\ph_r$ into a power series
$F(w_1,\ldots,w_r)$, defined as
\beglab{eqdefsubstF}
F(\ti\ph_1,\ldots,\ti\ph_r) \defeq \sum_{k\in\N^r} 
c_k \, \ti\ph_1^{k_1} \cdots \ti\ph_r^{k_r}
\edla
for $F = \sum_{k\in\N^r} c_k \, w_1^{k_1} \cdots w_r^{k_r}$.
One of the main results of \cite{KS} is
%
%
\begin{thm}    \label{thmsubstgR}
Let $r\ge1$ be an integer.
Then, for any convergent power series $F(w_1,\ldots,w_r) \in \C\{w_1,\ldots,w_r\}$
and for any resurgent series $\ti\ph_1,\ldots,\ti\ph_r \in \gR$ without constant term,
\[
F(\ti\ph_1,\ldots,\ti\ph_r) \in \gR.
\]
\end{thm}
It is the aim of this announcement to outline the idea of the proof,
based on the notion of ``discrete filtered sets'' and quantitative
estimates for the convolution of endless continuable functions.
%

Recall that the convolution in $\C[[\ze]]$, defined as
$\hat\ph_1 * \hat\ph_2 \defeq 
\cB\big( \cB\ii(\hat\ph_1) \cdot \cB\ii(\hat\ph_2) \big)$ 
(\ie the mere counterpart via~$\cB$ of multiplication of formal series without
constant term), takes the following form for convergent power series:
\begin{multline*}
\hat\ph_1,\hat\ph_2 \in \C\{\ze\} \quad\Rightarrow\quad
\hat\ph_1*\hat\ph_2(\ze) = \int_0^\ze \hat\ph_1(\ze_1)
\hat\ph_2(\ze-\ze_1)\,\dd\ze_1 \\
\ens\text{for}\;
\abs{\ze} < \min\{ \RCV(\hat\ph_1), \RCV(\hat\ph_2) \},
\end{multline*}
where $\RCV(\,\cdot\,)$ is a notation for the radius of convergence of
a power series.
The symbol~$\de$ which appears in~\eqref{eqdefgR} and~\eqref{eqdefcB}
is nothing but the convolution unit (obtained from $(\C[[\ze]],*)$ by adjunction of unit);
convolution makes $\C\de\oplus\C[[\ze]]$ a ring (isomorphic to the
ring of formal series $\C[[z]]$), of which $\C\de\oplus\C\{\ze\}$ is a subring.


It is proved in \cite{DO} that the convolution product of two
endlessly continuable functions is endlessly continuable, hence 
$\C\de\oplus\hat\gR$ is a subring of $\C\de\oplus\C\{\ze\}$
and $\gR$ is a subring of $\C[[z]]$.
However, to reach the conclusions of Theorem~\ref{thmsubstgR}, precise
estimates on the convolution product of an arbitrary number of
endlessly continuable functions are needed,
so as to prove the convergence of the series of holomorphic functions 
$\sum c_k \, \hat\ph_1^{*k_1} * \cdots * \hat\ph_r^{*k_r}$
and to check its endless continuability.


\begin{rem}    \label{remOmcontrestr}
Let $\Om$ be a nonempty closed discrete subset of~$\C$.
If a convergent power series~$\hat\ph$ meets the requirement of
Definition~\ref{defendlesscont} with 
\beglab{eqdefOmLrestr}
F_L = \{\, \om \in \Om \mid \abs{\om} \le L \,\}
\quad\text{for each $L>0$,}
\edla 
then it is said to be \emph{$\Om$-continuable}.
Pulling back this definition by~$\cB$, we obtain the definition of
an \emph{$\Om$-resurgent series}, which is a particular case of
resurgent series.
It is proved in \cite{SFunkEkva} that the space of $\Om$-continuable
functions is closed under convolution if and only if $\Om$ is stable
under addition.
Under that assumption, the space of $\Om$-resurgent series is thus a
subring of $\C[[z]]$; nonlinear analysis with $\Om$-resurgent series
is dealt with in \cite{NLresur}, where an analogue of
Theorem~\ref{thmsubstgR} is proved for them.
\end{rem}


Our study of endlessly continuable functions and our proof of
Theorem~\ref{thmsubstgR} are based on the notion of
``$\Om$-continuable function'', where $\Om$ is a ``discrete
  filtered set'' in the sense of \cite{DO};
this is a generalization of the situation described in
Remark~\ref{remOmcontrestr}, so that the meaning
of $\Om$-continuability
will now be more extended ($\Om$ will stand for a family of finite
sets $(\Om_L)_{L\ge0}$ not necessarily of the form~\eqref{eqdefOmLrestr}).

The plan of the paper is as follows.
Discrete filtered sets and the corresponding $\Om$-continuable
functions are defined in Section~\ref{sec:dfs}, 
where we give a refined version of the main result, Theorem~\ref{thmsubstOmgR}.
Then, in Section~\ref{sec:pfthm}, we state Theorem~\ref{thm:1.12}
which gives precise estimates for the convolution product of an
arbitrary number of $\Om$-continuable functions, and show how this
implies Theorem~\ref{thmsubstOmgR}.
Finally, in Section~\ref{sec:isot}, we sketch the main step of
the proof of Theorem~\ref{thm:1.12}.

\section{Discrete filtered sets and $\Om$-resurgent series}\label{sec:dfs}

\subsection{$\Om$-continuable functions and $\Om$-resurgent series}


We use the notations 
\[ 
\N = \{0,1,2,\ldots\}, \qquad
\Rp = \{\la\in\R\ |\ \la\geq0\}.
\]
\begin{dfn}[Adapted from \cite{CNP},\cite{DO}]
A \emph{discrete filtered set}, or \emph{\dfs}\ for short, is a family $\Om = (\Om_L)_{L\in\Rp}$ of
subsets of~$\C$ such that
\begin{enumerate}[i)]
\item
%
$\Om_L$ is a finite set for each $L\in\Rp$, 
\item
$\Om_{L_1}\subseteq \Om_{L_2}$ for $L_1\leq L_2$,
\item
there exists $\de>0$ such that $\Om_\de=\O$.
\end{enumerate}
Given a \dfs~$\Om$, we set 
\[
\cS_\Om \defeq \big\{ (\la,\om)\in\Rp\times\C \mid \om\in\Om_\la \big\}
\]
and define $\ov\cS_\Om$ as the closure of $\cS_\Om$ in 
$\Rp\times\C$. We then call
\[
\cM_\Om \defeq \big(\Rp\times\C\big) \setminus \ov\cS_\Om
\quad \text{(open subset of $\Rp\times\C$)}
\]
the \emph{allowed open set} associated with~$\Om$.
\end{dfn}
%
%
\begin{dfn}   \label{defOmallowedpaths}
  We denote by~$\Pi$ the set of all Lipschitz paths
  $\ga \col [0,t_*] \to \C$ such that $\ga(0)=0$, with some real $t_*\ge0$ depending on~$\ga$.
We then denote by $\ga_{|t} \in \Pi$ the restriction of~$\ga$ to the
interval $[0,t]$ for any $t\in[0,t_*]$, and by $L(\ga)$ the total
length of~$\ga$.
Given a \dfs~$\Om$, we call \emph{$\Om$-allowed path} any $\ga\in\Pi$ such
that
\[
\big( L(\ga_{|t}), \ga(t) \big) \in \cM_\Om
\ens\text{for all $t$.}
\]
We denote by $\Pi_\Om$ the set of all $\Om$-allowed paths.
\end{dfn}
%
%
\begin{dfn}   \label{defOmcontOmres}
Given a \dfs~$\Om$, we call \emph{$\Om$-continuable function} a
holomorphic germ $\hat\ph \in \C\{\ze\}$ which can be analytically
continued along any path $\ga\in\Pi_{\Om}$.
We denote by~$\hat\gR_\Om$ the set of all $\Om$-continuable functions
and define
\[
\gR_\Om \defeq \cB\ii \big( \C\de \oplus \hat\gR_\Om \big)
\]
to be the set of \emph{$\Om$-resurgent} series.
\end{dfn}
%

\begin{exa}   \label{exacloseddiscOm}
Given a closed discrete subset~$\Om$ of~$\C$, the formula
$\ti\Om_L \defeq \{\, \om \in \Om \mid \abs{\om} \le L \,\}$ for
$L\in\Rp$ defines a \dfs~$\ti\Om$.
Then the notion of $\Om$-continuability defined in
Remark~\ref{remOmcontrestr} agrees with the notion of
$\ti\Om$-continuability of Definition~\ref{defOmcontOmres}.
Thus we can identify the set~$\Om$ and the \dfs~$\ti\Om$.
\end{exa}


The relation between the $\Om$-continuable functions or the
$\Om$-resurgent series of Definition~\ref{defOmcontOmres} and the
resurgent series or the endlessly continuable functions of
Section~\ref{sec:0} is as follows:


\begin{prp}   \label{propidgROmgR}
A formal series $\ti\ph\in\C[[z]]$ is resurgent if and only if there
exists a \dfs~$\Om$ such that $\ti\ph$ is $\Om$-resurgent.
Correspondingly,
\beglab{eqidentitygROmgR}
\hat\gR = \bigcup_{\Om\in\{\text{\dfs}\}} \hat\gR_\Om.
\edla
\end{prp}


The proof of Proposition~\ref{propidgROmgR} is outlined in Section~\ref{sec:pfprpRes}.


\begin{rem}    \label{reminclusdfs}
Observe that
\[
\Om \subset \Om' \quad\Rightarrow\quad
\gR_\Om \subset \gR_{\Om'},
\]
where the symbol $\subset$ in the left-hand side stands for the
partial order defined by $\Om_L \subset \Om'_L$ for each~$L$.
Indeed, $\Om\subset\Om'$ implies $\cS_\Om\subset\cS_{\Om'}$, hence
$\cM_{\Om'} \subset \cM_\Om$ and $\Pi_{\Om'} \subset \Pi_\Om$.
\end{rem}


\begin{rem}    
Obviously, entire functions are always $\Om$-continuable:
$\gO(\C) \subset \hat\gR_\Om$ for all \dfs~$\Om$
(e.g.\ because $\gO(\C) = \hat\gR_{\O}$, denoting by~$\O$ the trivial
\dfs).
%
The inclusion is not necessarily strict for a non-trivial \dfs\ 
In fact, one can show that
\[
\hat\gR_\Om = \gO(\C)
\quad \Leftrightarrow \quad
\forall L>0,\; \exists L'>L\; \text{such that} \;
\Om_{L'} \subset \{\, \om\in\C \mid \abs{\om} < L \,\}.
\]
%
%
A simple example where this happens is when $\Om_L=\O$ for $0\le L<2$
and $\Om_L = \{1\}$ for $L\ge2$.
\end{rem}


\subsection{Sums of discrete filtered sets}   \label{sec:sumsdsf}


The proof of the following lemma is easy and left to the reader.

\begin{lmm}
Let~$\Om$ and~$\Om'$ be two \dfs\ 
Then the formula
\[
(\Om *\Om')_L \defeq \{\, \om_1+\om_2 \mid
\om_1\in\Om_{L_1}, \om_2\in\Om'_{L_2}, L_1+L_2=L \, \}
\cup\Om_{L}\cup\Om'_{L}
\quad\text{for $L\in\Rp$}
\]
defines a \dfs\ $\Om*\Om'$.
The law~$*$ on the set of all \dfs\ is commutative and associative.
The formula $\Om^{*n} \defeq 
\underbrace{ \Om* \cdots*\Om }_{\text{$n$ times}}$
(for $n\ge1$) defines an inductive
system and
\[
\Om^{*\infty} \defeq \varinjlim_n\ \Om^{*n}
\]
is a \dfs
\end{lmm}


We call $\Om*\Om'$ the \emph{sum} of the \dfs~$\Om$ and~$\Om'$.
In \cite{DO} the following is claimed:


\begin{thm}[\cite{DO}]\label{thm:DO}
  Assume that $\Om$ and $\Om'$ are \dfs\ and $\ti\ph\in\gR_\Om$,
  $\ti\psi\in\gR_{\Om'}$.
Then the product series $\ti\ph \cdot\ti\psi$ is $\Om*\Om'$-resurgent.
\end{thm}


In view of Proposition~\ref{propidgROmgR}, a direct consequence of
Theorem~\ref{thm:DO} is


\begin{crl}
  The space of resurgent formal series~$\gR$ is a subring of the ring
  of formal series $\C[[z]]$.
\end{crl}


Similarly, in view of Proposition~\ref{propidgROmgR}, Theorem~\ref{thmsubstgR}
is a direct consequence of
\begin{thm}   \label{thmsubstOmgR}
Let $r\ge1$ be integer and let $\Om_1$, \ldots, $\Om_r$ be \dfs\ 
Then for any convergent power series $F(w_1,\ldots,w_r) \in
\C\{w_1,\ldots,w_r\}$
and for any $\ti\ph_1,\ldots,\ti\ph_r \in \C[[z]]$ without constant
term, one has 
\[
\ti\ph_1 \in \gR_{\Om_1}, \ldots, \ti\ph_r \in \gR_{\Om_r}
\quad \Rightarrow \quad
F(\ti\ph_1,\ldots,\ti\ph_r) \in \gR_{\Om^*},
\]
where $\Om^* \defeq (\Om_1 * \cdots * \Om_r)^{*\infty}$.
\end{thm}


The proof of Theorem~\ref{thmsubstOmgR} is outlined in
Sections~\ref{sec:pfthm} and~\ref{sec:isot}.

Theorem~\ref{thm:DO} may be viewed as a particular case of
Theorem~\ref{thmsubstOmgR} (by taking $F(w_1,w_2)=w_1 w_2$).
The proof of the former theorem consists in checking that, for
$\hat\ph\in\hat\gR_\Om$ and $\hat\psi \in \hat\gR_{\Om'}$, the
convolution product $\hat\ph*\hat\psi$ can be analytically continued
along the paths of $\Pi_{\Om*\Om'}$ and thus belongs to $\hat\gR_{\Om*\Om'}$.
In the situation of Theorem~\ref{thmsubstOmgR}, with the
notation~\eqref{eqdefsubstF}, we have 
$\hat\psi_k \defeq c_k \, \hat\ph_1^{*k_1} * \cdots * \hat\ph_r^{*k_r} \in \hat\gR_{\Om^*}$
for each nonzero $k\in\N^r$, but some analysis is required to prove the
convergence of the series $\sum\hat\psi_k$ of $\Om^*$-continuable
functions in~$\hat\gR_{\Om^*}$;
what we need is a precise estimate for the convolution product of
an arbitrary number of endlessly continuable functions, 
and this will be the content of Theorem~\ref{thm:1.12}.


In the particular case of a closed discrete subset~$\Om$ of~$\C$
assumed to be stable under the addition and viewed as a \dfs\ as in
Example~\ref{exacloseddiscOm}, we have $\Om^{*\infty} = \Om$
and Theorem~\ref{thm:1.12} was proved for that case in \cite{NLresur}.
One of the main purposes of \cite{KS}
is to extend the techniques of \cite{NLresur} to the more general
setting of endlessly continuable functions.


\subsection{Upper closure of a \dfs\ and sketch of proof of
  Proposition~\ref{propidgROmgR}} \label{sec:pfprpRes}


The proof of Proposition~\ref{propidgROmgR} makes use of the notion of
``upper closure'' of a \dfs, which allows to simplify a bit the
definition of $\Om$-allowedness for a path, and thus of
$\Om$-continuability for a holomorphic germ.

\begin{dfn}
We call \emph{upper closure} of a \dfs~$\Om$ the family of sets
$\ti\Om = (\ti\Om)_{L\ge0}$ defined by
\[
\ti\Om_L \defeq \bigcap_{\eps>0} \Om_{L+\eps}
\quad\text{for every $L\in\Rp$.}
\]
\end{dfn}

Notice that $\Om \subset \ti\Om$.


\begin{lmm}   \label{lemclosOm}
Let $\Om$ be a \dfs\ 
Then its upper closure~$\ti\Om$ is a \dfs, and there exists a real sequence
$(L_n)_{n\ge0}$ such that
$0=L_0 < L_1 < L_2 < \cdots$
and
\[
L_n < L < L_{n+1} 
\quad \Rightarrow \quad
\ti\Om_{L_n} = \ti\Om_L = \Om_L
\]
for every integer $n\ge0$.
\end{lmm}


\begin{lmm}   \label{lemOmallowedness}
Let $\Om$ be a \dfs. Then
\beglab{eqovcSOmcStiOm}
\ov\cS_\Om = \cS_{\ti\Om}.
\edla
Consequently, $\Om$-allowedness admits the following characterization:
for a path $\ga\in\Pi$,
\begin{align*}
%
%
\ga \in \Pi_\Om
\quad&\Leftrightarrow\quad
\text{for all $t$,}\;\, 
\ga(t) \in \C \setminus \ti\Om_{L(\ga_{|t})} \\[1ex]
&\Leftrightarrow\quad
\text{for all $t$,\; $\exists n$ such that}\;\,
L(\ga_{|t}) < L_{n+1} 
\ens\text{and}\ens
\ga(t) \in \C \setminus \ti\Om_{L_n}
\end{align*}
(with the notation of Lemma~\ref{lemclosOm}).
\end{lmm}


The proofs of Lemma~\ref{lemclosOm} and Lemma~\ref{lemOmallowedness}
are easy; see \cite{KS} for the details.


Notice that, given~$\Om$ and $\ga\in\Pi_\Om$, it may be impossible to
find \emph{one} real $L\ge L(\ga)$ such that $\ga(t) \in \C\setminus\ti\Om_L$
for \emph{all} $t>0$.
Take for instance~$\Om$ defined by $\Om_L \defeq \O$ for $0\le L < 2$
and $\Om_L \defeq \{1,2\}$ for $L\ge2$ (so $\Om = \ti\Om$ in
that case),  
and $\ga\in\Pi_\Om$ following the line segment $[0,3/2]$, then winding once
around~$2$ and ending at~$3/2$:
no ``uniform''~$L$ can be found for that path.
Observe that, in that example, there exists $\hat\ph\in\hat\gR_\Om$ with
an analytic continuation along~$\ga$ such that the resulting
holomorphic germ at~$3/2$ is singular at~$1$ (but~$1$ is not singular
for the principal branch of~$\hat\ph$).
This is why the proof of Proposition~\ref{propidgROmgR} requires a bit
of work.


\begin{proof}[Sketch of the proof of Proposition~\ref{propidgROmgR}]
It is sufficient to prove~\eqref{eqidentitygROmgR}.
%
%
Suppose $\hat\ph \in \hat\gR_\Om$ for a certain \dfs~$\Om$ and let
$L>0$.
Then~$\hat\ph$ meets the requirement of
Definition~\ref{defendlesscont} with $F_L = \ti\Om_L$,
hence $\hat\ph\in\hat\gR$.
Thus $\hat\gR_\Om \subset \hat\gR$.
%

Suppose now $\hat\ph\in\hat\gR$. 
In view of Definition~\ref{defendlesscont}, we have $\de \defeq
\RCV(\hat\ph) >0$ and, for each positive integer~$n$, we can choose a
finite set~$F_{n}$ such that
\begin{eq-text}   \label{eqpropertyFnde}
the germ~$\hat\ph$ can be analytically continued along any path
$\ga\col[0,1]\to\C$ of~$\Pi$ such that $L(\ga) < (n+1)\de$ and
$\ga\big( (0,1] \big) \subset \C\setminus F_{n}$.
\end{eq-text}
Let $F_0\defeq \O$. The property~\eqref{eqpropertyFnde} holds for
$n=0$ too.
For every real $L\ge0$, we set
\[
\Om_L \defeq \bigcup_{k=0}^n F_{k}
\qquad \text{with $n \defeq \floor{L/\de}$.}
\]
One can check that $\Om\defeq (\Om_L)_{L\in\Rp}$ is a \dfs\ which coincides with its upper
closure.
In \cite{KS}, it is shown with the help of Lemma~\ref{lemOmallowedness} that $\hat\ph\in\hat\gR_\Om$.
\end{proof}

\section{The Riemann surface~$X_\Om$ -- an estimate
  for the convolution product of several $\Om$-continuable functions} 
\label{sec:pfthm}

As announced in Section~\ref{sec:sumsdsf}, we will now state a theorem
from which Theorem~\ref{thmsubstOmgR} and thus
Theorem~\ref{thmsubstgR} follow.
A few preliminaries are necessary.


In all this section we suppose that~$\Om$ is a fixed \dfs\


\begin{dfn}[
Adapted from \cite{DO}]
We call \textit{$\Om$-endless Rieman surface} any triple $(X,\fp,\uO)$ such
that
$X$ is a connected Riemann surface,
$\fp \col X \to \C$ is a local biholomorphism,
$\uO \in \fp\ii(0)$,
and any path $\ga \col [0,t_*] \to \C$ of~$\Pi_\Om$ has a lift
$\uga \col [0,t_*] \to X$ such that $\uga(0) = \uO$.
\end{dfn}

Notice that, given $\ga \in \Pi_\Om$, the lift~$\uga$ is unique
(because the fibres of~$\fp$ are discrete).


It is shown in \cite{KS} how, among all $\Om$-endless Riemann
surfaces, one can construct an object $(X_{\Om},\fp_\Om,\uO_\Om)$ which is cofinal in the
following sense:

\begin{prp}   \label{lemuniversalXOm}
There exists an $\Om$-endless Riemann surface $(X_{\Om},\fp_\Om,\uO_\Om)$ such
that, for any $\Om$-endless Riemann surface $(\ti X,\ti\fp,\ti\uO)$,
there is a local biholomorphism $\fq \col X_\Om \to \ti X$ such
that $\fq(\uO_\Om)=\ti\uO$ and the diagram
\[
\begin{xy}
(0,20) *{X_{\Om}},
(30,20) *{\ti X},
(15,0) *{\C},
{(6,20) \ar (24,20)},
{(2,16) \ar (13,4)},
{(28,16) \ar (17,4)},
(14,23) *{\fq},
(3,10) *{\fp_\Om},
(27,10) *{\ti\fp}
\end{xy}
\]
is commutative.
The $\Om$-endless Riemann surface $(X_{\Om},\fp_\Om,\uO_\Om)$ is unique up to
isomorphism
and $X_\Om$ is simply connected.
\end{prp}

The reader is referred to \cite{KS} for the proof.
Notice that in particular, for any \dfs~$\Om'$ such that
$\Om\subset\Om'$ as in Remark~\ref{reminclusdfs},
%
Proposition~\ref{lemuniversalXOm} yields a local biholomorphism $\fq \col X_{\Om'}\to X_{\Om}$.


For an arbitrary connected Riemann surface~$X$, we denote by~$\gO_X$
the sheaf of holomorphic functions on~$X$.
If $\fp \col X \to \C$ is a local biholomorphism, then there is a natural morphism
$\fp^* \col \fp\ii\gO_{\C} \to \gO_X$.
Recall that $\fp\ii\gO_{\C}$ is a sheaf on~$X$, whose stalk at a
point $\un\ze_* \in \fp\ii(\ze_*)$ is
%
$\big( \fp\ii\gO_{\C} \big)_{\un\ze_*} = \gO_{\C,\fp(\un\ze_*)} \cong \C\{\ze-\ze_*\}$.


\begin{lmm}
Let $\hat\ph \in \C\{\ze\} = \gO_{\C,0}$. 
Then the following properties are equivalent:
\begin{enumerate}[{\rm i)}]
\item
$\hat\ph\in\gO_{\C,0}$ is $\Om$-continuable,
\item
$\fp_\Om^*\hat\ph \in \gO_{X_{\Om},\uO_\Om}$ can be analytically continued
along any path on $X_{\Om}$,
\item
$\fp_\Om^*\hat\ph \in \gO_{X_{\Om},\uO_\Om}$ can be extended to 
$\Ga(X_{\Om},\gO_{X_{\Om}})$.
\end{enumerate}
\end{lmm}

So, the morphism~$\fp_\Om^*$ allows us to identify an $\Om$-continuable
function with a function holomorphic on the whole of the Riemann
surface~$X_\Om$:
\[
\restr{\fp_\Om^*}{\hat\gR_\Om} \col
\hat\gR_\Om \isom \Ga(X_{\Om},\gO_{X_{\Om}}).
\]
The reader is referred to \cite{KS} for the details.


\begin{nota}
For $\de>0$ small enough so that $\Om_\de=\O$ and for $L>0$, we set
\begin{gather*}
\Pi_\Om^{\de,L} \defeq \big\{\, \ga\in\Pi \mid
L(\ga) \le L \;\;\text{and}\;\;
\abs{ L(\ga_{|t})-\la }^2 + \abs{ \ga(t)-\om }^2 \ge \de^2
\;\, \text{for all $(\la,\om)\in\cS_\Om$ and $t$} \,\big\}, \\[1ex]
K_{\Om}^{\de,L} \defeq \big\{\, \un\ze\in X_\Om \mid 
\exists t_*\ge0 \;\,\text{and}\;\, \exists \ga \col [0,t_*] \to \C 
\;\text{path of $\Pi_\Om^{\de,L}$ such that 
$\un\ze = \uga(t_*)$}
\,\big\}.
\end{gather*}
\end{nota}


The condition $\Om_\de=\O$ is meant to ensure that
$\Pi_\Om^{\de,L}$ is a nonempty subset of $\Pi_\Om$ and hence
$K_\Om^{\de,L}$ is a nonempty compact subset of~$X_\Om$.
In fact, with the notation of Lemma~\ref{lemclosOm}, $\Om_\de=\O$ as
soon as $\de<L_1$ and then, for every $\ga\in\Pi$,
\[
L(\ga)+\de \le L_1 
\quad \Rightarrow \quad
\ga \in \Pi_\Om^{\de,L} 
\ens \text{for all $L \ge L(\ga)$.}
\]
(In particular, if $L+\de\le L_1$, then $\Pi_\Om^{\de,L}$ consists
exactly of the paths of~$\Pi$ which have length $\le L$.)
Notice that, if $\Om_\de = \O$, then $\Om^{*n}_\de=\O$ for all $n\ge1$
and $\Om^{*\infty}_\de=\O$.

On the other hand, for a given \dfs~$\Om$, the family
$(K_\Om^{\de,L})_{\de,L>0}$ yields an exhaustion of~$X_\Om$ by compact
subsets.


We are now ready to state the main result of this section,
which is the analytical core of our study of the convolution of
endelssly continuable functions:

\begin{thm}\label{thm:1.12}
Let $\Om$ be a \dfs\ and let $\de,L>0$ be reals such that $\Om_{2\de}=\O$.
Then there exist $c,\de'>0$ such that $\de'\le\de$ and, for every integer $n\ge1$ and for
every
$\hat f_1, \ldots, \hat f_n \in \hat\gR_\Om$,
the function $1*\hat f_{1}*\cdots *\hat f_{n}$ belongs to $\hat\gR_{\Om^{*n}}$
and satisfies the following estimates:
\beglab{1.11}
\sup_{K_{\Om^{*n}}^{\de,L}} \abs*{
\fp_{\Om^{*n}}^* \big( 1*\hat f_{1}*\cdots *\hat f_{n}\big )
}
\leq
\frac{c ^n}{n!}
\ \sup_{\substack{ L\pp1,\ldots,L\pp n>0 \\ L\pp 1+\cdots+L\pp n=L}}
\ \sup_{K_{\Om}^{\de',L\pp 1}} \abs{ \fp_\Om^{*} \hat f_1 }
\, \cdots
\sup_{K_{\Om}^{\de',L\pp n}} \abs{ \fp_\Om^{*} \hat f_n }.
\edla
\end{thm}


The main step of the proof of Theorem~\ref{thm:1.12} is sketched in
Section~\ref{sec:isot} 
(possible values for~$c$ and~$\de'$ are indicated in~\eqref{eqdefdepcthmestim};
the point is that they depend on~$\Om$, $\de$ and~$L$, but not on~$n$
nor on $\hat f_1, \ldots, \hat f_n$).
The full proof of Theorem~\ref{thm:1.12} is in \cite{KS}.


\begin{proof}[Theorem~\ref{thm:1.12} implies
  Theorem~\ref{thmsubstOmgR}]
Let $r\ge1$ and let $\Om_1,\ldots,\Om_r$ be \dfs\
Let 
$\Om \defeq \Om_1 * \cdots * \Om_r$,
so that $\Om_i \subset \Om$ and $\hat\gR_{\Om_i} \subset \hat\gR_\Om$ for
$i=1,\ldots,r$.

Let $F \in \C\{w_1,\ldots,w_r\}$. Denote its coefficients by
$(c_k)_{k\in\N^r}$ as in~\eqref{eqdefsubstF} and pick $C,\La>0$ such
that
$\abs{c_k} \le C \La^{\abs{k}}$ for all~$k$,
with the notation $\abs{k} \defeq k_1+\cdots+k_r$.

Let $\ti\ph_1,\ldots\ti\ph_r$ be formal series without constant term
such that $\hat\ph_i \defeq \cB\ti\ph_i \in \gR_{\Om_i}$ for each~$i$. 
We get the uniform convergence of
$\sum c_k\, \fp_{\Om^{*\infty}}^* (1*\hat\ph_1^{*k_1}*\cdots*\hat\ph^{*k_r})$
on every compact subset of~$X_{\Om^{*\infty}}$ by means of estimates
of the form
\[
\sup_{K_{\Om^{*\infty}}^{\de,L}} \abs{ c_k\,\fp_{\Om^{*\infty}}^* (1*\hat\ph_1^{*k_1}*\cdots*\hat\ph^{*k_r}) }
\le
\sup_{K_{\Om^{*n}}^{\de,L}} \abs{ c_k\,\fp_{\Om^{*n}}^* (1*\hat\ph_1^{*k_1}*\cdots*\hat\ph^{*k_r}) }
\le
C \La^n \cdot \frac{c^n}{n!} \cdot M^n,
\]
where $\de,L>0$ with $\Om_{2\de}=\O$, $n\defeq \abs{k} \ge1$, 
$M \defeq \sup_{1\le i\le r}\, \sup_{K_{\Om}^{\de',L}} \abs{ \fp_\Om^{*} \hat\ph_i }$,
and the positive reals~$c$ and~$\de'$ stem from Theorem~\ref{thm:1.12}.
\end{proof}

\section{Sketch of the proof of Theorem~\ref{thm:1.12} -- construction
  of adapted deformations of the identity}
\label{sec:isot}
%
\subsection{Preliminaries}
%

Let $\Om$ be a \dfs\ 
%
%
Recall from Definition~\ref{defOmallowedpaths}
and~\ref{defOmcontOmres} that $\Om$-continuability is defined by means
of the allowed open subset $\cM_\Om$ of $\Rp\times\C$ associated
with~$\Om$, and of paths $\ga\in\Pi$ which satisfy
\[
\ti\ga(t) \defeq \big( L(\ga_{|t}), \ga(t) \big) \in \cM_\Om
\quad\text{for all $t$.}
\]
The latter condition is $\ga\in\Pi_\Om$.
Conversely, notice that
\begin{eq-text}   \label{eqtextcharacOmalltiga}
if $t \in [0,t_*] \mapsto 
\ti\ga(t) = \big( \la(t),\ga(t) \big) \in \cM_\Om$
is a piecewise $C^1$ path such that $\ti\ga(0)=(0,0)$ and 
$\la'(t) = \abs{\ga'(t)}$ for a.e.~$t$, then $\ga\in\Pi_\Om$.
\end{eq-text}
Recall from Section~\ref{sec:pfthm} that, if $\ga \col [0,t_*] \to \C$
is a path of~$\Pi_\Om$, then there is a unique path
$\uga \col [0,t_*] \to X_\Om$ such that $\uga(0)=\uO_\Om$ and
$\fp_\Om\circ\uga = \ga$.

We fix $\rho>0$ such that $\Om_{2\rho}=\O$ and set
\[ U \defeq \{\, \ze\in\C \mid \abs{\ze} < 2\rho \,\}. \]
For every $\ze\in U$, the path $\ga_\ze \col t\in[0,1] \mapsto t\ze$
is in~$\Pi_\Om$ and the formula
\[
\ze \in U \mapsto \gL(\ze) \defeq \uga_\ze(1) \in X_\Om
\]
defines a holomorphic section~$\gL$ of~$\fp_\Om$ on~$U$.
Let $\uU \defeq \gL(U)$: this is an open subset of~$X_\Om$
containing~$\uO_\Om$ and we have mutually inverse biholomorphisms
\[
\gL \col U \isom \uU, \qquad 
\restr{\fp_\Om}{\uU} \col \uU \isom U.
\]


\begin{nota}	\label{notesimplexn}
For any $n\ge1$, we denote by~$\De_n$ the $n$-dimensional
simplex
\[
\De_n \defeq \{\, (s_1,\ldots,s_n)\in\R^n \mid s_1,\ldots,s_n\ge0
\;\text{and}\;
s_1 +\cdots + s_n \le 1 \,\}
\]
with the standard orientation, and by $[\De_n]\in\gE_n(\R^n)$ the corresponding integration
current.
%
%
For every $\ze\in U$, we consider the map
\[
\gD(\ze) \col
\vec s = (s_1,\ldots,s_n) \mapsto 
\gD(\ze,\vec s) \defeq \big( \gL(s_1\ze),\ldots,\gL(s_n\ze) \big) \in
\uU^n \subset X_\Om^n,
\]
defined in a neighbourhood of~$\De_n$ in~$\R^n$,
and denote by $\gD(\ze)_\# [\De_n] \in \gE_n(X_\Om^n)$ the push-forward
of~$[\De_n]$ by~$\gD(\ze)$.
%
\end{nota}


The reader is referred to \cite{NLresur} for the notations and notions
related to integration currents.
Let $n\ge1$, $\hat f_1,\ldots, \hat f_n \in \hat\gR_\Om$ and 
$\hat g \defeq 1*\hat f_1 * \cdots * \hat f_n$. 
Our starting point is
\begin{lmm}[\cite{NLresur}]
Let
\[ \be \defeq (\fp_\Om^*\hat f_1)\big(\un\ze_1\big) \cdots (\fp_\Om^*\hat f_n)\big(\un\ze_n\big) \,
\dd\un\ze_1 \wedge \cdots \wedge \dd\un\ze_n, \]
where we denote by $\dd\un\ze_1 \wedge \cdots \wedge \dd\un\ze_n$ the pullback
of the $n$-form $\dd\ze_1 \wedge \cdots \wedge \dd\ze_n$
in $X_\Om^n$ by $\fp^{\otimes n} \col X_\Om^n \to \C^n$.
Then 
\[
\ze \in U \quad\Rightarrow\quad
\hat g(\ze) =
\gD(\ze)_\# [\De_n](\be).
\]
\end{lmm}

%
\subsection{$\ga$-adapted deformations of the identity}
%

Let $L>0$ and $\de \in (0,\rho)$. Let $\ga \col [0,1] \to \C$ be a path of
$\Pi_{\Om^{*n}}^{\de,L} \subset \Pi_{\Om^{*n}}$.
We want to study the analytic continuation of~$\hat g$ along~$\ga$,
which amounts to studying 
$(\fp_{\Om^{*n}}^*\hat g)\big( \uga(t) \big)$,
where~$\uga$ is the lift of~$\ga$ in $X_{\Om^{*n}}$ which starts at~$\uO_{\Om^{*n}}$.
Without loss of generality, we may assume
that there exists $a\in(0,1)$ such that 
\[
0 < \abs{\ga(a)} < \rho, \qquad
\ga(t) = \tfrac{t}{a}\ga(a) \ens\text{for $t\in[0,a]$,}
\qquad \text{$\restr{\ga}{[a,1]}$ is $C^1$.}
\]


\begin{dfn}
For $\ze \in \C$ and $1\le i\le n$, we set
\begin{align*}
\cN(\ze) &\defeq \big\{\, \big(\un\ze_1,\ldots,\un\ze_n\big) \in X_\Om^n \mid 
\fp_\Om\big(\un\ze_1\big) + \cdots + \fp_\Om\big(\un\ze_n\big) = \ze \,\big\},
\\[1ex]
\cN_i &\defeq \big\{\, \big(\un\ze_1,\ldots,\un\ze_n\big) \in X_\Om^n \mid \un\ze_i = \uO_\Om \,\big\}.
\end{align*}
We call \emph{$\ga$-adapted deformation of the identity on~$\un V$}
any family $(\Psi_t)_{t\in[a,1]}$ of maps 
\[
\Psi_t \col \un V \to X_\Om^n,
\]
where~$\un V$ is a neighbourhood of $\gD\big( \ga(a) \big)(\De_n)$ in~$X_\Om^n$,
such that 
$\Psi_a = \ID$, 
the map
$\big(t,\vuze\big) \in [a,1] \times \un V \mapsto
\Psi_t\big(\vuze\big) \in X_\Om^n$
is locally Lipschitz,
and for any $t\in [a,1]$ and $i=1,\ldots,n$,
\begin{align*}
\vuze \in \cN\big(\ga(a)\big) \quad & \Rightarrow \quad
\Psi_t\big(\vuze\big) \in \cN\big(\ga(t)\big), \\
\vuze \in \cN_i \quad & \Rightarrow \quad
\Psi_t\big(\vuze\big) \in \cN_i.
\end{align*}
\end{dfn}

The above notion is a slight generalization of the ``$\ga$-adapted
origin-fixing isotopies'' which appear in \cite[Def.~5.1]{NLresur}.
Adapting the proof of \cite[Prop.~5.2]{NLresur}, we get
\begin{prp}[\cite{NLresur}]
If $(\Psi_t)_{t\in[a,1]}$ is a $\ga$-adapted deformation of the
identity, then
\beglab{eqconthatgugat}
(\fp_{\Om^{*n}}^*\hat g)\big( \uga(t) \big) = 
\big( \Psi_t \circ \gD\big( \ga(a) \big) \big)_\# [\De_n](\be)
\qquad\text{for $t\in[a,1]$.}
\edla
\end{prp}
Notice that, with the notations 
\beglab{eqdefunzeitzeit}
\big(
\un\ze_1^t, \ldots, \un\ze_n^t
\big) \defeq
\Psi_t\circ \gD\big( \ga(a) \big) \col \De_n \to X_\Om^n,
\qquad \ze_i^t \defeq \fp_\Om \circ \un\ze_i^t \quad\text{for $1\le
  i\le n$,}
\edla
formula~\eqref{eqconthatgugat} can be rewritten as 
\beglab{eqconthatgugatbis}
(\fp_{\Om^{*n}}^*\hat g)\big( \uga(t) \big) = 
\int_{\De_n} 
(\fp_\Om^*\hat f_1)\big(\un\ze_1^t\big) \cdots (\fp_\Om^*\hat f_n)\big(\un\ze_n^t\big)
\Det\bigg[ \frac{\pa\ze_i^t}{\pa s_j}\bigg]_{1\le i,j\le n}
\, \dd s_1 \cdots \dd s_n
\edla
(for each~$t$, the partial derivatives $\frac{\pa\ze_i^t}{\pa s_j}$
exist almost everywhere on~$\De_n$ by Rademacher's theorem, for the functions $\vec s \mapsto
\ze_i^t(\vec s)$ are Lipschitz).

%
\subsection{Sketch of the proof of Theorem~\ref{thm:1.12}}
%

%
We define a function $\eta\ge0$ by the formula
\[
v = (\la,\xi) \in \cM_\Om \mapsto 
\eta(v) \defeq \dist\big( (\la,\xi), \{ (0,0) \} \cup \ov\cS_\Om \big),
\]
where $\dist(\cdot,\cdot)$ is a notation for the Euclidean distance in $\R\times\C
\simeq \R^3$.

The following three lemmas allow to prove Theorem~\ref{thm:1.12}.
The reader is referred to \cite{KS} for their proofs.


\begin{lmm}
The function~$D$ defined by the formula
\[
D\big(t,(v_1,\ldots,v_n)\big) \defeq
\eta(v_1) + \cdots + \eta(v_n) + 
\dist\big( \big( L(\ga_{|t}),\ga(t) \big), v_1+\cdots+v_n \big)
\]
is everywhere positive on $[a,1] \times \cM_\Om^n$ and the formula
\beglab{eqdefvecX}
\vec X(t,\vec v) =  \left| \begin{aligned}
X_1 &\defeq \frac{\eta(v_1)}{D(t,\vec v)} \big( \abs{\ga'(t)}, \ga'(t) \big)\\
& \qquad \vdots \\[1ex]
X_n &\defeq \frac{\eta(v_1)}{D(t,\vec v)} \big( \abs{\ga'(t)}, \ga'(t) \big)
\end{aligned} \right.
\edla
defines a non-autonomous vector field 
$\vec X(t,\vec v) \in T_{\vec v} \big( \cM_\Om^n \big) \simeq (\R\times\C)^n$
on $[a,1]\times\cM_\Om^n$, 
which admits a flow map $\Phi_t \col \cM_\Om^n \to \cM_\Om^n$ between time~$a$ and time~$t$
for every $t\in [a,1]$.
\end{lmm}


\begin{lmm}
One can define a $\ga$-adapted deformation of the identity $(\Psi_t)_{t\in[a,1]}$ on~$\uU^n$
as follows:
for every $\vuze = \big( \gL(\ze_1),\ldots \gL(\ze_n) \big) \in \uU^n$
and $i \in \{1,\ldots,n\}$,
we set $v_j \defeq (\abs{\ze_j},\ze_j)$ for each $j\in\{1,\ldots,n\}$
and define a path $\ti\ga_{i,\vuze} \col [0,1] \to \cM_\Om$ by
\[
t \in [0,a] \ens\Rightarrow\ens 
\ti\ga_{i,\vuze}(t) \defeq \tfrac{t}{a} v_i,
\qquad
t \in [a,1] \ens\Rightarrow\ens 
\ti\ga_{i,\vuze}(t) \defeq \pi_i \circ \Phi_t(v_1,\ldots, v_n),
\]
where $\pi_i \col \cM_\Om^n \to \cM_\Om$ is the projection onto the
$i^{\text{th}}$ factor;
then, by virtue of~\eqref{eqtextcharacOmalltiga}, 
the $\C$-projection of $\ti\ga_{i,\vuze}$ is 
a path $\ga_{i,\vuze} \in \Pi_\Om$, and we set
\[
\Psi_t\big(\vuze\big) \defeq \big(
\uga_{1,\vuze}(t), \ldots, \uga_{n,\vuze}(t)
\big) \in X_\Om^n
\quad \text{for $t\in[a,1]$.}
\]
\end{lmm}


\begin{lmm}   \label{lemestimdet}
Consider the functions
$\vec s \in \De_n \mapsto
v_i^a(\vec s) \defeq \big(
s_i\abs{\ga(a)}, s_i\ga(a)
\big) \in \cM_\Om$ for $1\le i \le n$ and, 
for each $t\in [a,1]$,
\[
(v_1^t, \ldots, v_n^t) \defeq \Phi_t \circ (v_1^a, \ldots, v_n^a)
\col \De_n \to \cM_\Om^n.
\]
Suppose $|\ga(a)|>\de$ and let
\[
\cM_\Om^{\de,L} \defeq \big\{\,
(\la,\ze) \in \Rp\times\C \mid
\la\leq L \;\;\text{and}\;\;
\dist\big( (\la,\ze), u \big) \ge \de
\;\, \text{for all $u\in\{ (0,0) \} \cup \ov\cS_\Om$}
\,\big\}.
\]
Then
\begla
(v_1^t, \ldots, v_n^t)(\De_n) \subset
\bigcup_{L\pp 1+\cdots+L\pp n=L(\ga_{|t})}
\cM_\Om^{\de(t),L\pp 1}
\times\cdots\times
\cM_\Om^{\de(t),L\pp n},
\edla
where
$\de(t) \defeq 
%
\frac12 \rho\,\ee^{-2\sqrt{2}\de\ii(L(\ga_{|t}) - \abs{\ga(a)})}$.
Moreover, for each~$i$, the $\C$-projection of~$v_i^t$ is a Lipschitz
function $\ze_i^t \col \De_n \to \C$ and the almost everywhere defined
partial derivatives 
$\frac{\pa\ze_i^t}{\pa s_j}$
satisfy
\[
\abs*{\Det\bigg[ \frac{\pa\ze_i^t}{\pa s_j}\bigg]_{1\le i,j\le n}} 
\le \big( c(t) \big)^n,
\]
where
$c(t) \defeq 
\abs{\ga(a)} \,\ee^{3\sqrt{2}\de\ii(L(\ga_{|t}) - \abs{\ga(a)})}$.
\end{lmm}


We set
\beglab{eqdefdepcthmestim}
\de' \defeq 
%
\frac12 \rho\,\ee^{-2\sqrt{2}\de\ii L},
\qquad
c \defeq 
\abs{\ga(a)} \,\ee^{3\sqrt{2}\de\ii L},
\edla
so that $\de' \le \de(t)$
and $c \ge c(t)$ for all $t\in[a,1]$.
Theorem~\ref{thm:1.12} follows from the previous estimates
and the identity~\eqref{eqconthatgugatbis}
(the functions~$\ze_i^t$ in~\eqref{eqdefunzeitzeit} and in
Lemma~\ref{lemestimdet} are indeed the same).



\newpage
\vspace{.4cm}

\noindent {\em Acknowledgements.}
{This work has been supported by Grant-in-Aid for JSPS Fellows Grant Number 15J06019, French National Research Agency reference ANR-12-BS01-0017 and Laboratoire Hypathie A*Midex.}



\vspace{2.3cm}

\noindent
Shingo \textsc{Kamimoto}\\
Department of Mathematics,
Hiroshima University, Hiroshima 739-8526, Japan.

\vspace{1cm}

\noindent
David \textsc{Sauzin}\\
Laboratorio Fibonacci, CNRS--CRM Ennio De
  Giorgi SNS Pisa, Italy.\\
{\itshape Current address}: IMCCE, CNRS--Observatoire de Paris, France.

\vspace{2.3cm}

\begin{flushright}

\textit{\small August 21, 2015. Revised February 25, 2016.}

\end{flushright}

\end{document}